\documentclass{amsart}

\usepackage{amsmath}
\usepackage{amssymb}
\usepackage{amscd}
\usepackage{amsthm}
\usepackage{enumerate}
\usepackage{array}
\usepackage{epsf} 
\usepackage{pictexwd}
\usepackage{rotating}
\usepackage{float}  \restylefloat{figure} 
\usepackage[OT2,OT1]{fontenc}

\theoremstyle{plain}
\newtheorem*{thmA}{Theorem A}
\newtheorem*{thmB}{Theorem B}
\newtheorem{lem}{Lemma}
\newtheorem{prop}[lem]{Proposition}

\theoremstyle{definition}

\theoremstyle{remark}
\newtheorem*{rem}{Remark}

\newdimen\XX

\newcommand{\noline}{\\[-\extrarowheight]}

\def\su{\mathfrak{su}}

\def\forcehmode{\hskip0pt\relax}

\let\myskip=\medskip

\font\msbm=msbm10
\def\semiprod{\hbox{\msbm\char111}}

\def\st{\,\,\big|\,\,}
\def\bs{\backslash}
\def\<{\langle}
\def\>{\rangle}
\def\ie{i.e.\ }

\let\ge=\geqslant
\let\le=\leqslant
\let\emptyset=\varnothing

\def\definebb#1=#2.{\def#1{{{\mathbb #2}^{\vphantom{x}}}}}
\definebb \c=C.
\definebb \r=R.
\definebb \d=D.
\definebb \n=N.
\definebb \z=Z.
\definebb \q=Q.

\def\calF{\mathcal F}

\DeclareMathOperator{\Cl}{Cl}

\DeclareMathOperator{\Id}{Id}
\DeclareMathOperator{\Int}{Int}
\DeclareMathOperator{\Isom}{Isom}
\DeclareMathOperator{\PSL}{PSL}
\DeclareMathOperator{\SO}{SO}
\DeclareMathOperator{\MathOpPSU}{PSU}
\DeclareMathOperator{\MathOpSU}{SU}

\let\Re=\undefined \DeclareMathOperator{\Re}{Re}

\def\dd{\partial}
\def\ddP{{\dd P}}

\def\al{{\alpha}}
\def\ve{{\varepsilon}}
\def\thet{{\vartheta}}
\def\De{{\Delta}}

\def\PSU{\MathOpPSU(1,1)}
\def\SU{\MathOpSU(1,1)}
\def\tSU{\widetilde{\MathOpSU}(1,1)}

\def\G{{\Gamma}}
\def\bG{\overline\G}
\def\la{{\lambda}}

\def\Grp{G}
\def\tGrp{\tilde\Grp}
\def\Lrp{L}
\def\tLrp{\tilde\Lrp}

\def\bdlmap{s}
\def\covbdlmap{s}
\def\covmap{\pi}

\def\WZZW{\begin{pmatrix} w&z \\ \bar z&\bar w \end{pmatrix}}

\def\Eg{E_g}
\def\Ig{I_g}
\def\Hg{H_g}

\def\Fg{F_g}
\def\Eh{E_h}
\def\Fh{F_h}
\def\Ia{I_a}
\def\Ee{E_e}
\def\Ie{I_e}
\def\He{H_e}
\def\Xe{X_e}
\def\Fe{F_e}

\def\calFg{\calF_g}

\def\barg{{\bar g}}
\def\bare{{\bar e}}

\def\bE{{\bar E}}
\def\bI{{\bar I}}

\def\bEbg{\bE_\barg}
\def\bIbg{\bI_\barg}
\def\bEbe{\bE_\bare}
\def\bIbe{\bI_\bare}

\def\Tx{T(x)}
\def\Tu{T(u)}
\def\Qx{Q_x}
\def\Qu{Q_u}
\def\Pu{P_u}
\def\Xu{X_u}

\def\deck{r_d}

\def\spieg{\eta}
\def\inv{\ve}

\def\Rx{R_x}

\def\Gu{{\G(u)}}
\def\Guou{{\Gu\backslash\{u\}}}

\def\cupddQxGu{{\bigcup\limits_{\hbox to 0pt{\hss$\scriptstyle x\in\Gu$\hss}}\dd\Qx}}
\def\capRxGu{{\bigcap\limits_{\hbox to 0pt{\hss$\scriptstyle x\in\Gu$\hss}}\Rx}}
\def\cupQxGu{{\bigcup\limits_{\hbox to 0pt{\hss$\scriptstyle x\in\Gu$\hss}}\Qx}}
\def\capRxGuou{{\bigcap\limits_{\hbox to 0pt{\hss$\scriptstyle x\in\Guou$\hss}}\Rx}}
\def\cupIntQxGuou{{\bigcup\limits_{\hbox to 0pt{\hss$\scriptstyle x\in\Guou$\hss}}\Int\Qx}}

\def\hF{\hat F}

\newenvironment{pictexure}[2]{\bgroup
  \def\PICdirectory{#1}\beginpicture
  \input #1/#2.pic \ignorespaces}{\endpicture\egroup}

\newcommand{\Pictexure}[2]{\begin{pictexure}{#1}{#2}\relax\end{pictexure}}

\let\cdlabelsize=\small
\def\cdlhss{0.5em\relax}
\newdimen\defaultarrow \defaultarrow 28pt
\newdimen\arrskip \arrskip=0.11\defaultarrow
\newif\ifcdarrowdepth \cdarrowdepthtrue
\newif\ifcdarrowheight \cdarrowheighttrue
\def\spacedfill #1#2{
  \hskip#1\arrskip\hbox to#1\defaultarrow{\cdlabelsize#2}\hskip#1\arrskip}

\def\cdarrow#1#2{
   \ifnum#1=0  \spacedfill{#2}{\rightarrowfill}\else
   \ifnum#1=180\spacedfill{#2}{\leftarrowfill }\else
   $\vcenter{\hbox{\vtop{\hrule height0pt
  \hbox{\begin{turn}{#1}%
   \hbox{\spacedfill{#2}{\rightarrowfill}}%
  \end{turn}}}}}$\fi\fi}

\def\cdarrowxx#1#2#3#4{
  \putcdlabel {\cdlabelsize$\displaystyle\genfrac{}{}{0pt}{}{#1}{#2}$}
   onto box {\cdarrow{#3}{#4}}}

\def\cdarrowv#1#2#3#4{
  \llap{\cdlabelsize$#1$\hskip\cdlhss}%
  \cdarrow{#3}{#4}%
  \rlap{\cdlabelsize\hskip\cdlhss$#2$}}

\newbox\arrowbox
\def\putcdlabel #1 onto box #2{
  \setbox\arrowbox=\hbox{#2}%
  \setbox0=\hbox{\rlap{\copy\arrowbox}\hbox to\wd\arrowbox{\hss #1\hss}}%
  \ifcdarrowheight\ht0=\ht\arrowbox\relax\fi \cdarrowheighttrue
  \ifcdarrowdepth\dp0=\dp\arrowbox\relax\fi \cdarrowdepthtrue
  \box0}

\begin{document}

\author[Anna Pratoussevitch]{Anna Pratoussevitch}
\address{Mathematisches Institut\\ Universit\"at Bonn\\ Beringstra{\ss}e~1 \\ 53115 Bonn}
\email{anna@math.uni-bonn.de}

\title{Fundamental Domains in Lorentzian Geometry}

\begin{date}  {\today} \end{date}

\thanks{Research partially supported by the
Graduiertenkolleg Mathematik in Bonn, financed by DFG}

\begin{abstract} 
We consider discrete subgroups $\Gamma$ of the simply connected Lie
group $\widetilde{\operatorname{SU}}(1,1)$ of finite level, 
i.e.\ the subgroup intersects the centre of $\widetilde{\operatorname{SU}}(1,1)$
in a subgroup of finite index, this index is called the level of the group. 
The Killing form induces a Lorentzian metric of constant curvature 
on the Lie group $\widetilde{\operatorname{SU}}(1,1)$.
The discrete subgroup $\Gamma$ acts on $\widetilde{\operatorname{SU}}(1,1)$ by left translations. 
We describe the Lorentz space form
$\widetilde{\operatorname{SU}}(1,1)/\Gamma$ by
constructing a fundamental domain $F$ for $\Gamma$. 
We want $F$ to be a polyhedron with totally geodesic faces. 
We construct such $F$ for all $\Gamma$ satisfying the following condition: 
The image $\bar\Gamma$ of $\Gamma$ in $\operatorname{PSU}(1,1)$
has a fixed point $u$ in the unit disk of order larger than the index of
$\Gamma$. 
The construction depends on the group $\Gamma$ and on the orbit $\Gamma(u)$ of
the fixed point~$u$.

 \end{abstract}









\maketitle

\section{Introduction}

We consider the universal cover of $\PSU\cong\PSL(2,\r)$,
the group of orien\-ta\-tion-preserving isometries of the hyperbolic plane.
Here our model of the hyperbolic plane is the unit disk $\d$ in $\c$.

\myskip
The kernel of the universal covering map $\tSU\to\PSU$
is the centre $Z$ of the group $\tSU$,
an infinite cyclic group.
Therefore, for each natural number $k$ there is a unique connected $k$-fold
covering of $\PSU$.
For $k=2$ this is the group
$$\SU=\left\{\WZZW\st(w,z)\in\c^2,~|w|^2-|z|^2=1\right\}.$$

\myskip
The level of a discrete subgroup $\G\subset\tSU$
is the index of $\G\cap Z$ as a subgroup of $Z$.
There is a one-to-one correspondence between discrete subgroups of level~$k$ in $\tSU$ 
and liftings of discrete subgroups in $\PSU$ into the $k$-fold covering of $\PSU$.

\myskip
We consider a discrete subgroup $\G$ in $\tSU$ of finite level $k$.
We suppose that the image $\bG$ of $\G$ in $\PSU$ has 
at least one fixed point in $\d$ of order $p$,
\ie a point in $\d$, which is fixed by a nontrivial element of $\bG$ of order $p$.
Furthermore we assume that the order $p$ of the fixed point is larger then the
level $k$ of the subgroup $\G$.
Our construction depends on the choice of the fixed point $u\in\d$ of $\bG$,
or actually on its orbit $\bG(u)$.

\myskip
The Killing form on the Lie group $\tSU$ gives rise
to a Lorentz biinvariant metric of constant curvature.
The quotient of $\tSU$ by the discrete subgroup $\G$ is a Lorentz space form
with respect to this metric,
\ie a complete Lorentz manifold of constant curvature
(compare R.S.~Kulkarni and F.~Raymond~\cite{KR}).

\myskip
The main result of this paper is the construction of fundamental domains
for the action of $\G$ on $\tSU$ by left translations,
applicable to any discrete subgroup $\G$ in $\tSU$ as above.
This fundamental domain is a polyhedron in the
Lorentz manifold $\tSU$ with totally geodesic faces.
For a co-compact subgroup the corresponding fundamental domain is compact.
The precise formulation of this result is contained in Theorems A and B.

\myskip
The construction of fundamental domain was a part of author's Ph.D. thesis~\cite{Pr:2001}.
It was outlined in the survey article~\cite{BPR}.
The object of this paper is to provide complete proofs for this result
together with a concise and self-contained description of the construction.

\myskip
Our results generalize a construction by Th.~Fischer~\cite{Fi}.
He suggested how to construct a fundamental domain for the action of
a discrete subgroup of $\PSU$ by left multiplication.
His construction can be interpreted in our terms as a construction for
discrete subgroup of $\tSU$ of level~$1$.
A less technical proof of Th.~Fischer's result is given in~\cite{Ba}. 

\myskip
The study of discrete subgroups of finite level was originally motivated
by some deep connections between these subgroups and quasi-homogeneous
isolated singularities of complex surfaces studied by J.~Milnor, I.~Dolgachev and W.~Neumann~\cite{Mi, Do83, Neu77, Neu83}.
In particular the quotient $\tSU/\G$ is diffeomorphic to the link of some
quasi-homogeneous Gorenstein singularity.
For a more detailed treatment of this connection see \cite{BPR}, \S 1--2.

\myskip
The paper is organized as follows:
We start in section~\ref{analogons} by discussion of low-dimensional analogues of our problem 
and use these examples to indicate some of the main ingredients of our construction.
Section~\ref{prelim} contains some general remarks on the Lie groups $\SU$ and $\tSU$
and their embeddings in some $4$-dimensional pseudo-Euclidean space resp.\ in a certain $\r_+$-bundle, 
the universal cover of a positive cone in that pseudo-Riemannian space.
We describe in section~\ref{elements} some elements of the construction, such as
affine half-spaces and their substitutes in the $\r_+$-bundle.
We also define prismatic sets $\Qu$, certain finite intersections of half-spaces, 
and study their properties. 

\myskip
After that we are prepared to state in section~\ref{results} our main results,
Theorems A and B, and to prove them.
In section~\ref{examples} we report on our explicit computations of fundamental
domains for certain infinite series of discrete subgroups.
The choice of this series is motivated by the connections between them and
some series of quasi-homogeneous surface singularities.
We also give some pictures of fundamental domains.
Finally, in section~\ref{outlook} we discuss some relations to and similarities with
other fundamental domain constructions and give an outlook on possible generalizations.

\myskip
The results described in section~\ref{examples} have been announced 
in~\cite{BPR} together with computations of fundamental domains 
for other singularities by E.~Brieskorn and F.~Rothenh\"ausler.
The synopsis of these results points at a very regular pattern for these
series of singularities and their fundamental domains.

\myskip
I am indebted to Egbert Brieskorn and Thomas Fischer who opened up this field of research.
I am very grateful to Egbert Brieskorn for his guidance and helpful discussions.
I would like to thank
Ludwig Balke,
Werner Ballmann,
Ilya Dogolazky,
Pierre Pansu,
and Frank Rothenh\"auser
for useful conversations related to this work.
I thank Ilya Dogolazky for his help in producing the figures.

\section{Low-Dimensional Analogues}

\label{analogons}

\noindent
Before we start to describe the construction for $\tSU$,
we discuss in this section the corresponding problem of finding fundamental domains for
an action of a discrete subgroup for two toy cases, 
for the one-dimensional Lie groups $\SO(2)$ and $\SO(1,1)$.
Some of the main ideas of the construction for $\tSU$ can be seen already in
the discussion of these low-dimensional examples.

\myskip
We first consider the Lie group $\SO(2)$.
We identify $\SO(2)$ with the circle
$$G=\{z(t)=(\cos t,\sin t)\st t\in\r\}$$
in the Euclidean plane $E^2.$
We consider a discrete subgroup $\G_m$, 
the finite cyclic subgroup of order $m$ generated by $z(2\pi/m)$,
and its action on $\SO(2)$ by left multiplication,
which extends to an action on $E^2$ by isometries.

\myskip
Clearly, the segment of length $2\pi/m$ with midpoint $z(0)$ is a fundamental domain for this action,
this is the Dirichlet domain with respect to the point $z(0)$.
However, for the description of Dirichlet domains we need the distance.
The following description of the same fundamental domain as projection of an affine
construction with tangent half-planes is more appropriate for generalizations
in pseudo-Riemannian setting.

\begin{figure}
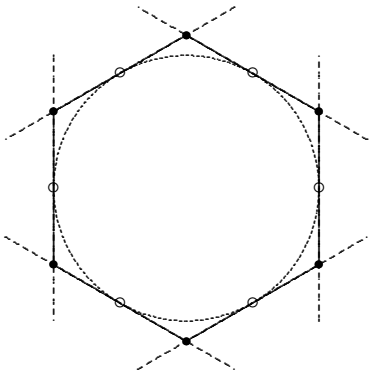

  \begin{center}
    \forcehmode
      \bgroup
        \beginpicture
          \input so2.tex
        \endpicture
      \egroup
  \end{center}
  \caption {Fundamental domain construction in $\SO(2)$}
  \label{so2-figa}
\end{figure}

\myskip
For $g\in G$ let $\Hg=\{a\in E^2\st\<a,g\>\le1\}$ be the half-plane 
with boundary tangent to the circle $G$ in the point $g$.
Then the intersection of the tangent half-planes
$$P=\bigcap\limits_{g\in\G_m}\Hg$$
is a regular $\G_m$-invariant $m$-gon,
its faces are fundamental domains for the action of $\G_m$ on the boundary $\dd P$,
and the projection of the faces under the contraction $a\mapsto a/|a|$
yields to a tiling of the circle by (Dirichlet) fundamental domains with respect to~$\G_m$.
Figure~\ref{so2-figa} illustrates the construction for $m=6$.


\myskip
A more involved example, where the construction of Dirichlet fundamental domains can be described in the same way,
is the case of the action of discrete subgroups on $\MathOpSU(2)$,
in particular the construction of the $4$-dimensional regular polyhedron
bounded by 120 dodecahedra, the tiling of the $3$-dimensional sphere
by 120 spherical dodecahedra and the resulting construction
of the Poincar\'e ho\-mo\-lo\-gy sphere using the binary icosahedral group.
Here we identify $\MathOpSU(2)$ with the $3$-dimensional sphere in the Euclidean space $E^4$.

\myskip
Our second description of the fundamental domain in $\SO(2)$ 
does not use any distance. 
It only uses the embedding of $\SO(2)$ in the (pseudo-)Euclidean space.
However, the simple-minded attempt to generalize this affine construction 
to pseudo-Riemannian quadrics fails.
Our second one-dimensional example, the hyperbola in the Minkowski plane,
shows, why the naive approach fails and what can be done about this.

\myskip
We now consider the Lie group $\SO(1,1)$ and identify it with the hyperbola
$$G=\{z(t,\ve)=\ve\cdot(\sinh t,\cosh t)\st t\in\r,~\ve=\pm1\}$$
in the Minkowski plane $E^{1,1}$ with metric induced by $\<a,b\>=a_1b_1-a_2b_2$.
Let us fix $d>0$ and consider a discrete subgroup $\G_d$, 
the subgroup generated by the elements $z(d,1)$ and $z(d,-1)$ and isomorphic to $\z\times\{\pm1\}$.
Moreover we consider the action of $\G_d$ on $\SO(1,1)$ by left multiplication,
which extends to an action on $E^{1,1}$ by isometries.
Clearly, the segment of length $d$ with midpoint $z(0,+1)$ is a fundamental domain for this action.
This is the Dirichlet domain with respect to the point $z(0,+1)$.


\begin{figure}
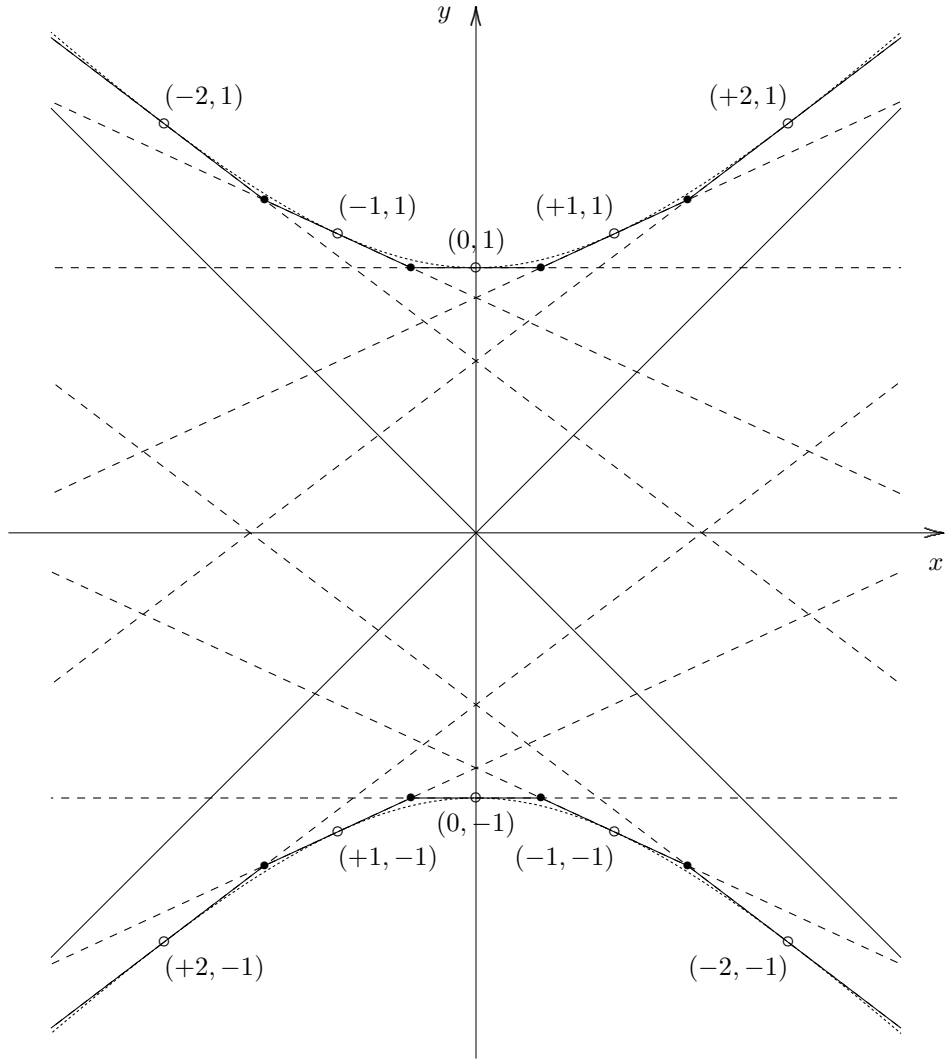

  \begin{center}
    \forcehmode
    \Pictexure{perl}{so11}
  \end{center}
  \caption {Fundamental domain construction in $\SO(1,1)$}
  \label{so11-figa}
\end{figure}

\myskip
Let $L$ be the cone over $G$, \ie $L=\r_+\cdot G=\{a\in E^{1,1}\st\<a,a\><0\}$.
For $g\in G$ let $\Hg$ be the half-plane $\Hg=\{a\in E^{1,1}\st\<a,g\>\ge-1\}$  
with boundary tangent to the hyperbola $G$ in the point $g$.
The polytope $P=\cap_{g\in\G_d}\Hg$ contains in this case only one point $(0,0)$,
so we can not get any information about the tiling from this set.
Instead we consider the set
$$P=\bigcup\limits_{k\in\z}\bigcap\limits_{\ve=\pm1} H_{z(k\cdot d,\ve)}$$
and the part of its boundary contained in $L$.

\myskip
Figure~\ref{so11-figa} shows the polyhedron $P$ and illustrates our
construction of fundamental domains for $\G_d$ on the part of the boundary of
$P$ lying over $G$, \ie in the positive cone $L$ over $G$.
The following statements can be easily verified:
The set $\dd P\cap L$ is $\G_d$-invariant,
and its faces are fundamental domains for the action of $\G_d$ on $\dd P\cap L$.
The projection $\dd P\cap L\to G$ given by the contraction $a\mapsto a/\sqrt{-(a,a)}$
is a $\G_d$-equivariant homeomorphism,
and the projection of the faces of $\dd P\cap L$ yields to a tiling of the
hyperbola by fundamental domains with respect to $\G_d$.

\myskip
This construction suggests some important ingredients of the construction for $\SU$,
namely the embedding of $\SU$ as a quadric
in a $4$-dimensional pseudo-Euclidean space,
appropriate decomposition of the discrete subgroup in countable many
finite subsets $\Tx$, $x\in X$,
and finally the study of the $4$-dimensional polyhedron
$$P=\bigcup\limits_{x\in X}\bigcap\limits_{g\in\Tx} \Hg.$$

\myskip
Some new ideas come in when we generalize the fundamental domain construction for $\tSU$.
We consider an embedding of $\tSU$ as an image of a section in a (trivial)
$\r_+$-bundle over $\tSU$, namely in the universal cover of the positive cone $\r_+\cdot\tSU$,
and we define appropriate substitutes for tangent spaces and half-spaces there.

\myskip
We also want to point out that these two one-dimensional examples are not only examples, 
we meet them again in the construction for $\SU$. 
We identify $\SU$ with a quadric in a $4$-dimensional pseudo-Euclidean space.
The subgroups $\SO(2)\cong\MathOpSU(1)$ and $\SO(1,1)$ of $\SU$
can be identified with certain plane sections of this quadric,
and the corresponding constructions of fundamental domains
are then sections of the construction for $\SU$.

\section{Preliminaries}

\label{prelim}

\noindent
We consider the $4$-dimensional pseudo-Euclidean space $E^{2,2}$ of signature $(2,2)$.
We think of $E^{2,2}$ as real vector space $\c^2\cong\r^4$ with the symmetric bilinear form
$$\<(z_1,w_1),(z_2,w_2)\>=\Re(z_1\bar z_2-w_1\bar w_2).$$

\myskip
In $E^{2,2}$ we consider the pseudo-hyperbolic space
\begin{align*}
 \Grp
  &=\left\{a\in E^{2,2}\st\<a,a\>=-1\right\}\\
  &=\left\{(z,w)\in E^{2,2}\st|z|^2-|w|^2=-1\right\}\\
\end{align*}
and the cone over $\Grp$
\begin{align*}
 \Lrp
  &=\r_+\cdot\Grp\\
  &=\left\{a\in E^{2,2}\st\<a,a\><0\right\}\\
  &=\left\{(z,w)\in E^{2,2}\st|z|<|w|\right\}.\\
\end{align*}
The bilinear form on $E^{2,2}$ induces a pseudo-Riemannian metric of signature~$(2,2)$
on $\Lrp$ and a Lorentz metric of signature $(2,1)$ on $\Grp$.

\myskip
We may think of $\Lrp$ as a $\r_+$-bundle over $\Grp$ with radial projection
$\bdlmap:\Lrp\to\Grp$ as bundle map.
The map $\Lrp\to\d$ defined by $(z,w)\mapsto z/w$ is principal $\c^*$-bundle, 
where the action of $\la\in\c^*$ is defined by $\la\cdot(z,w)=(\la^{-1}z,\la^{-1 }w)$.
Let $\covmap:\tGrp\to\Grp$ be the universal covering.
Henceforth we identify the Lie group $\SU$ with $\Grp$ via
$$\WZZW\mapsto(z,\bar w),$$
and $\tSU$ with $\tGrp$.
The biinvariant metrics on $\Grp$ and $\tGrp$
are proportional to the Killing forms.
We denote the pull-back $\tLrp\to\tGrp$
of the $\r_+$-bundle $\bdlmap:\Lrp\to\Grp$ under the covering map
$\covmap:\tGrp\to\Grp$ also by $\covbdlmap$.
The following diagram commutes
$$
  \begin{CD}
   \tLrp            @>\pi>> \Lrp          \\
   @V{\covbdlmap}VV         @VV{\bdlmap}V \\
   \tGrp            @>\pi>> \Grp          \\
  \end{CD}
$$
$\Grp$ resp.\ $\tGrp$ is canonically embedded in $\Lrp$ resp.\ $\tLrp$ and
therefore there exist canonical trivializations $\Lrp\cong\Grp\times\r_+$  resp.\ $\tLrp\cong\tGrp\times\r_+$.
The covering $\tLrp$ inherits canonically a pseudo-Riemannian metric from $\Lrp$.

\myskip
We now give a brief description of the full isometry group of $\tGrp$
(compare sections 2.1--2.3 in \cite{KR}).
The action of $\tGrp\times\tGrp$ on $\tGrp$ via
$$(g,h)\cdot x=gxh^{-1}$$
is by Lorentz isometries since the metric is biinvariant.
The standard estimates for the dimension of the isometry group
show that the identity component $\Isom_0(\tGrp)$ is isomorphic to $(\tGrp\times\tGrp)/\De_Z$, 
where $\De_Z=\{(z,z)\st z\in Z\}$ and $Z$ is the centre of $\tGrp$.
The full isometry group of $\tGrp$ has four components
corresponding to time- and/or space-reversals.
Let $\inv$ be the geodesic symmetry at the identity given by $g\mapsto g^{-1}$
and $\spieg$ the lift of the conjugation by the matrix
$\left(\begin{smallmatrix} 0&1 \\ 1&0 \end{smallmatrix}\right)$
in $\Grp$ fixing the identity.
Then $\inv$ preserves the space-orientation and reverses the time-orientation,
while $\spieg$ reverses both the space- and time-orientation.
Moreover, the group $\Isom^+(\tGrp)=\<\Isom_0(\tGrp),\spieg\>$
is the full group of orientation-preserving isometries and
$$
  \Isom(\tGrp)
  =\<\Isom_0(\tGrp),\spieg,\inv\>
  \cong\Isom_0(\tGrp)\semiprod(\<\spieg\>\times\<\inv\>)
$$
is the full isometry group of $\tGrp$.


\myskip
The universal coverings $\tGrp\subset\tLrp$ can also be described as follows:
\begin{align*}
  \tLrp&=\{(z,\al,r)\in\c\times\r\times\r_+\st|z|<r\} \\
  \intertext{and}
  \tGrp&=\{(z,\al,r)\in\tLrp\st |z|^2=r^2-1\}\approx\{(z,\al)\in\c\times\r\}.
\end{align*}
We call the number $\al\in\r$ the argument of the element $(z,\al,r)\in\tLrp$.
Then the maps $\covbdlmap:\tLrp\to\tGrp$ and $\covmap:\tLrp\to\Lrp$
can be described as
$$
  \covbdlmap(z,\al,r)
  =\left(\la^{-1}z,\al,\la^{-1}r\right)
  \quad\hbox{with}\quad\la=\sqrt{r^2-|z|^2}
$$
and
$$\covmap(z,\al,r)=(z,re^{i\al}).$$

\section{The Elements of the Construction}

\label{elements}

\noindent
For $g\in\tGrp$ let $\Eg$ resp.\ $\Ig$ be the connected component
of $\covmap^{-1}(\bEbg)$ resp.\ $\covmap^{-1}(\bIbg)$ containing~$g$, where
$\barg:=\covmap(g)$ is the image of $g$ in $\Grp$,
\begin{align*}
  \bEbg&:=\{a\in\Lrp\st\<g,a\>=-1\} \\
  \intertext{is the intersection of $\Lrp$ with the affine tangent space on $\Grp$ in the point $\barg$ and}
  \bIbg&:=\{a\in\Lrp\st\<g,a\>\le-1\}
\end{align*}
is the intersection of $\Lrp$ with the half-space of $\c^2$ bounded by $\bEbg$
and not containing~$0$.
$\bEbg$ and $\bIbg$ are simply connected and even contractible,
hence their pre-images under the covering map $\covmap$ consist of infinitely
many connected components, one of them containing~$g$.

\myskip
The three-dimensional submanifold $\Eg$ subdivides $\tLrp$ in two
connected components, the closure of one of them is $\Ig$, and we denote the
closure of the other by $\Hg$.
The boundary of $\Ig$ resp.~$\Hg$ is equal to $\Eg$.

\myskip
As an example, for the unit elements $e=(0,0,1)$ in~$\tGrp$ and $\bare=\covmap(e)=(0,1)$ in~$\Grp$,
we have
$$\bIbe=\{(z,w)\in\c^2\st\Re(w)\ge1,~|z|<|w|\},$$
the boundary $\bEbe$ of $\bIbe$ is a rotational hyperboloid of one sheet.
The pre-image of $\bIbe$ is
$$
  \covmap^{-1}(\bIbe)
  =\{(z,\al,r)\in\c\times\r\times\r_+\st r\cdot\cos\al\ge1,~|z|<r\}.
$$
The connected components of $\covmap^{-1}(\bIbe)$ resp.\ $\covmap^{-1}(\bEbe)$
containing~$e$ are
\begin{align*}
  \Ie
  &=\left\{
     (z,\al,r)\in\c\times\r\times\r_+
     \st
     |\al|<\frac{\pi}{2},~r\ge\frac{1}{\cos\al},~|z|<r
   \right\}\\
  \intertext{and}
  \Ee
  &=\left\{
     (z,\al,r)\in\c\times\r\times\r_+
     \st
     |\al|<\frac{\pi}{2},~r=\frac{1}{\cos\al},~|z|<r
   \right\}.
\end{align*}
The subsets $\Eg$ resp.\ $\Ig$ have the analogous properties
because $\Eg=g\cdot\Ee$ and $\Ig=g\cdot\Ie$.

\myskip
We make use of the following construction (compare \cite{Mi}).
Given a base-point $x\in\d$ and a real number $t$,
let $\rho_x(t)\in\PSU$ denote the rotation through angle $t$ about the point $x$.
Thus we obtain a homomorphism $\rho_x:\r\to\PSU$,
which clearly lifts to the unique homomorphism $r_x:\r\to\tSU$
into the universal covering group.
Since $\rho_x(2\pi)=\Id_{\d}$,
it follows that the lifted element $r_x(2\pi)$ belongs
to the central subgroup $Z$ of $\tSU$.
Note that this element $r_x(2\pi)\in Z$ depends continuously on $x$, and
therefore is independent of the choice of $x$.
We easily compute $r_0(2t)=(0,-t,1)$ 
and hence $r_x(2\pi)=r_0(2\pi)=(0,-\pi,1)$ for all $x\in\d$.
Moreover we obtain 
$$r_0(2t)\cdot(z,\al,r)=(ze^{it},\al-t,r).$$

\myskip
Let $\G$ be a discrete subgroup of finite level $k$ in $\tSU$
and let $\bG$ be the image of $\G$ in $\PSU$.
We assume the existence of a fixed point $u\in\d$ of $\bG$.
The isotropy group $\bG_u$ of $u$ in $\bG$ is a finite cyclic group
generated by $\rho_u(2\pi/p)$, where $p:=|\bG_u|$.
The isotropy group $\G_u$ of $u$ in $\G$ is a infinite cyclic group
generated by $\deck:=r_u(2\thet)$, where~$\thet:=\frac{\pi k}{p}$.
We can assume without lost of generality that $u=0\in\d$.
Under this assumption it follows
$$
  \deck=\left(0,-\thet,1\right)
  \quad\hbox{and}\quad
  \deck\cdot(z,\al,r)=\left(ze^{i\thet},\al-\thet,r\right).
$$
An important assumption for the following construction is
$$p>k,$$
the order of $u$ as a fixed point of~$\bG$ is greater than the level of the group.
In terms of $\deck$ this means that the argument $\thet$ of $\deck$ is less then $\pi$.

\myskip
Now let us start with the construction
of fundamental domains for the action of $\G$ on $\tGrp$.
For a point $x$ in the orbit $\G(u)$ let $\Tx$ be the left coset
$$\Tx=\{g\in\G\st g(u)=x\}$$ 
of the isotropy group $\G_u$ and let
$$\Qx=\bigcap\limits_{g\in\Tx}\Hg.$$

\begin{figure}
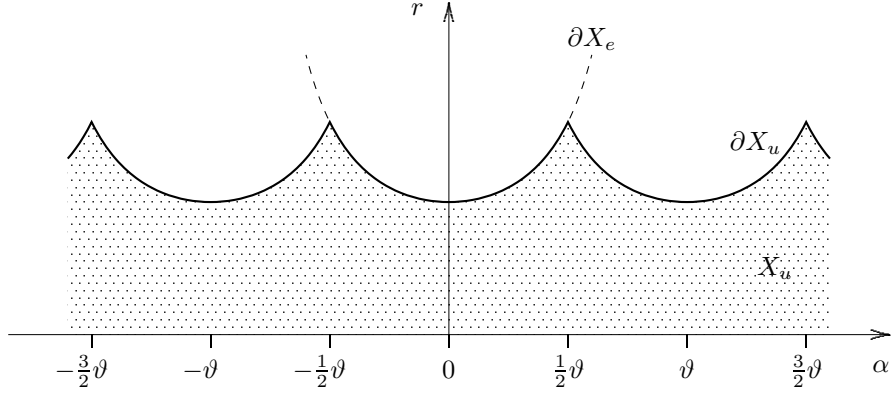

  \begin{center}
    \forcehmode
    \Pictexure{perl}{qu}
  \end{center}
  \caption {The image $X_u$ of $Q_u$ in the $(\alpha,r)$-half-plane}
  \label {figa}
\end{figure}

\myskip
As an example, for $x=u$ we have $\Tu=\G_u$,
the infinite cyclic subgroup generated by the element $\deck=(0,-\thet,1)$.
The generator $\deck$ acts on $\tGrp$ by left multiplication
$$\deck\cdot(z,\al,r) = (ze^{i\thet},\al-\thet,r)$$
and it acts on the $(\al,r)$-half-plane by the translation $\tau$ mapping $(\al,r)$ to $(\al-\thet,r)$.
How does the set $\Qu=\cap_{m\in\z}H_{\deck^m}$ look like?
The images of the sets $H_{\deck^m}=\deck^m\cdot\He$ under the projection $(z,\al,r)\mapsto(\al,r)$
are the translates $\tau^m(\Xe)$ of the image
$$\Xe=\{(\al,r)\in\r\times\r_+\st r\cdot\cos\al\le1\enskip\text{or}\enskip|\alpha|\ge{\pi}/{2}\}$$
of $\He$.
The manifold $\Qu$ is a disc bundle over its image
$\Xu=\bigcap_{m\in\z}\tau^m(\Xe)$ in the $(\al,r)$-plane.
The shaded area in figure~\ref{figa} is $\Xu$.


The manifolds $g\Qu$ play a central role in our construction.
We want to explain the geometric nature of these objects.
We have described $\Qu$ as a disc bundle over the set $\Xu$ in the $(\al,r)$-half-plane $\r\times\r_+$.
We may describe $\Qu\subset\tLrp\subset\c\times\r\times\r_+$ as
$$\Qu=(\c\times\Xu)\cap\tLrp.$$
We think of $\Xu$ as a universal covering of a punctured plane polygon.
Consider the following diagram of covering maps
$$
  \begin{tabular}{*2{c@{}}c}
    $\r\times\r_+$
          & \cdarrowheightfalse\cdarrowxx{\pi'}{}{0}{1}
          & $\c^\text{\rlap{$*$}}$ \\
    & \cdarrowxx{}{\pi~}{-45}{1.4142} & \!\cdarrowv{}{\pi''}{-90}{1} \\
    &                                   & $\c^\text{\rlap{$*$}}$
  \end{tabular}
$$
where $\pi(\al,r)=re^{i\al}$ and $\pi'(\al,r)=r^{1/k}e^{i\al/k}$ and $\pi''(z) = z^k$.
Consider the curve $\pi(\dd\Xu)$.
It is easy to see that this is a regular star polygon $\big\{\frac{2p}{k}\big\}$ when $k$ is odd
and a regular star polygon $\big\{\frac{p}{k}\big\}$ when $k$ is even.
(For the definition of a star polygon see for example H.S.M.~Coxeter~\cite{Co}, \S2.8, pp. 36--38.)
Therefore the curve $\pi'(\dd\Xu)$ is a curvilinear $2p$-gon
covering the star polygon once or twice.
Let $P'\subset\c$ and $P=\Pu\subset\c$ be the plane areas bounded by the curvilinear polygon
$\pi'(\dd\Xu)$ and by the star polygon $\pi(\Xu)$. The images of
$\Xu$ are the punctured plane polygons $\pi'(\Xu)=P'\setminus\{0\}$ and
$\pi(\Xu)=P\setminus\{0\}$.
We think of the product $\c\times P'$ as a 4-dimensional $2p$-gonal {\it prism}.
$\c\times\Xu$ is the universal covering of the pierced prism $\c\times(P'\setminus\{0\})$.
The product $\c\times P\subset\c^2$ might be considered as a 4-dimensional {\it ``star prism''.}
Its axis $\c\times\{0\}$ does not meet $\Lrp\subset\c\times\c^*$.
Therefore the universal covering $\pi:\tLrp\rightarrow\Lrp$ maps $\Qu$
to the intersection of $\Lrp$ with the star prism:
$$\pi(\Qu)=\Lrp\cap(\c\times\Pu).$$

\begin{rem}
With the help of the figure~\ref{figa} we easily prove the estimate
$$r\le\frac{1}{\cos\frac{\thet}{2}}$$
for points $(z,\al,r)\in\Qu$.
From this inequality we shall obtain in lemma~\ref{Qx}(i) some useful estimates
for points of $\Qx$.
\end{rem}


\myskip
The following properties can be easily proven for $x=u$
from the explicit description of the set $\Qu$.
The subsets $\Qx$ have the analogous properties
because $\Qx=a\cdot\Qu$ for $a\in\Tx$.

\begin{lem}
\label{Qx}
For a point $x\in\G(u)$ the following holds:
\begin{enumerate}[(i)]
\item
For any point $(z,w)\in\covmap(\Qx)$ we have
$$|w|-|z|\le|w-xz|\le f(|x|),$$
where
$$f(t):=\frac{\sqrt{1-t^2}}{\cos\frac{\thet}{2}}.$$
\item
The set $\Qx$ is a subgraph of a section
in the bundle $\tLrp\cong\tGrp\times\r_+$,
while its boundary is the graph of this section.
This means that for some section $s:\tGrp\to\r_+$ 
in the bundle $\tLrp$ we have
\begin{align*}
  \Qx&=\{(a,r)\in\tGrp\times\r_+\st r\le s(a)\},\\
  \dd\Qx&=\{(a,r)\in\tGrp\times\r_+\st r=s(a)\}
\end{align*}
\end{enumerate}
\end{lem}

\begin{lem}
\label{Qxlocfin}
The family $(\Qx)_{x\in\G(u)}$ is locally finite
in the sense that any point of $\tLrp$ has a neighbourhood intersecting
only finitely many prisms $\Qx$.
\end{lem}


\begin{proof}
We prove that the family $(\covmap(\Qx))_{x\in\G(u)}$ is locally finite (in
$\Lrp$). This fact implies the local finiteness of the family
$(\Qx)_{x\in\G(u)}$, since if a subset $U$ of $\Lrp$ has an empty intersection
with $\covmap(\Qx)$ then the intersection of the pre-image $\covmap^{-1}(U)$
with $\Qx$ is empty too.
By lemma~\ref{Qx}(i) for any point $x\in\G(u)$
and any point $(z,w)\in\covmap(\Qx)$ the difference $|w|-|z|$ is bounded from
above by $f(|x|)$.
The values $f(t)$ tend to zero as $t$ tends to $1$.
Choosing a point $(z_0,w_0)\in\Lrp$ and a positive number $\ve<|w_0|-|z_0|$,
the neighbourhood $U:=\{(w,z)\in\Lrp\st|w|-|z|>\ve\}$ of the point
$(z_0,w_0)$ can intersect $\covmap(\Qx)$ only for $|x|$ sufficiently
small (so that $f(|x|)>\ve$).
But the group $\G$ is discrete, so there are only finitely many points $x$ in
$\G(u)$ with norm $|x|$ under a given bound.
This finishes the proof.
\end{proof}

\begin{rem}
This property of $\Qx$ allows us to deal
with $P=\cup\Qx$ in a similar way as with a finite union of polytopes.
\end{rem}

\begin{lem}
\label{EgQxlocfin}
The family $(\Eg\cap Q_{g(u)})_{g\in\G}$ is locally finite.
\end{lem}

\begin{proof}
This is immediate from the local finiteness of the family $(\Qx)_{x\in\G(u)}$
plus the easy observation that the family $(\Eg\cap Q_{g(u)})_{g\in\G_u}$ is
locally finite.
\end{proof}

We consider in $\tLrp$ the four-dimensional polytope
$$
  P:=\bigcup\limits_{x\in\G(u)}\Qx
  =\bigcup\limits_{x\in\G(u)}\bigcap\limits_{g\in\Tx}\Hg.
$$

\begin{lem}
\label{ddPtGrp-homeo}
The projection~$\ddP\to\tGrp$ is a $\G$-equivariant homeomorphism.
\end{lem}


\begin{proof}
From lemma~\ref{Qx}(ii) we know that
the set $\Qx$ is a subgraph of a section
in the bundle $\tLrp\cong\tGrp\times\r_+$.
The polyhedron $P=\cup\Qx$ inherits this property of the prisms $\Qx$
as a union of a locally finite family of subgraphs.
But for a subgraph of a section in the bundle $\tLrp$
it is clear that the bundle map $\tLrp\to\tGrp$ induces a homeomorphism
from its boundary (equal to the graph of the section) onto $\tGrp$.
This homeomorphism is $\G$-equivariant
since the projection $\tLrp\to\tGrp$ is $\G$-equivariant.
\end{proof}

\section{The Main Results}

\label{results}

\noindent
Now we can state the main result

\begin{thmA}
The boundary of~$P$ is invariant with respect to the action of~$\G$.
The subset
$$\Fg=\Cl_{\ddP}(\Int(\dd\Hg\cap\ddP))$$
is a fundamental domain for the action of~$\G$ on~$\ddP$.
The family~$(\Fg)_{g\in\G}$ is locally finite in~$\ddP$.
The projection~$\tLrp\to\tGrp$ induces a $\G$-equivariant
homeomorphism~$\ddP\to\tGrp$.
The image~$\calFg$ of~$\Fg$ under the projection is a fundamental domain
for the action of~$\G$ on~$\tGrp$.
The family~$(\calFg)_{g\in\G}$ is locally finite.
For every pair of elements $g,h\in\G$ with~$g\ne h$ the intersection~$\calF_g\cap\calF_h$
lies in a totally geodesic submanifold of~$\tGrp$.
\end{thmA}

\begin{proof}
First let us prove that the intersection $\Int(\Fg)\cap\Fh$ is empty if $g\ne h$.
Suppose on the contrary that there are elements~$g,h\in\G$
such that $g\ne h$ and $\Int(\Fg)\cap\Fh\ne\emptyset$.
Since $\Eg\cap\dd P$ is closed, it follows by general topology arguments that
$\Int\Fg=\Int\Cl\Int(\Eg\cap\ddP)=\Int(\Eg\cap\ddP)$
(in this proof all closures are taken in~$\ddP$).
The assumption $\Int(\Fg)\cap\Fh\ne\emptyset$ can be rewritten as
$\Int(\Eg\cap\ddP)\cap\Cl\Int(\Eh\cap\ddP)\ne\emptyset$.
It can be shown in the standard way that if $A$ and~$B$ are closed subsets of
some topological space such that $\Int(A)\cap\Cl\Int(B)\ne\emptyset$
then it follows $\Int(A\cap B)\ne\emptyset$.
In our case this implies $\Int(\Eg\cap\Eh\cap\ddP)\ne\emptyset$.
But since the totally geodesic submanifolds~$\Eg$ and~$\Eh$ intersect transversally,
the intersection $\Eg\cap\Eh$ has no inner points in~$\ddP$.

Since $\Fg\subset\Eg\cap Q_{g(u)}$ lemma~\ref{EgQxlocfin} implies that
the family~$(\Fg)_{g\in\G}$ is locally finite in~$\ddP$.
Lemma~\ref{ddPtGrp-homeo} says that the projection~$\ddP\to\tGrp$ is a $\G$-equivariant homeomorphism.

Now let us prove that $\Cl(\cup_{g\in\G}\Int\Fg)=\ddP$.
Since
$$
  \Cl\big(\bigcup\limits_{g\in\G}\Int\Fg\big)
  \supset\bigcup\limits_{g\in\G}\Cl\Int\Fg
  =\bigcup\limits_{g\in\G}\Fg,
$$
it suffices to prove that $\cup_{g\in\G}\Fg=\ddP$.
Consider $a\in\ddP$.
From the definition of~$P$ and lemma~\ref{EgQxlocfin}
it follows that there exist elements $g_1,\dots,g_n\in\G$
and a neighbourhood~$U$ of the point~$a$ in~$\tLrp$ such that
$$\ddP\cap U=\bigcup\limits_{i=1}^n\,(E_{g_i}\cap\ddP\cap U).$$
We may assume without loss of generality that
the map $\covmap|_U:U\to\covmap(U)$ is a homeomorphism.
The image of $P\cap U$ under this homeomorphism is an intersection of an open
subset of~$\Lrp$ with a finite union of finite intersections of
half-spaces~$\Hg$ with the property $a\in\dd\Hg$.
Suppose that $a\not\in\Cl\Int(E_{g_i}\cap\ddP)=F_{g_i}$ for all $i\in\{1,\dots,n\}$.
This is only possible if for each~$i\in\{1,\dots,n\}$
the set $E_{g_i}\cap\ddP\cap U$ is contained in a $2$-dimensional submanifold of~$\tLrp$.
But from lemma~\ref{ddPtGrp-homeo} it follows that $\ddP$ is homeomorphic
to a $3$-dimensional manifold~$\tGrp$.
This contradiction implies that $a\in\Fg$ for some~$g\in\G$.
\end{proof}

\begin{rem}
The family $(\Qx)_{x\in\G(u)}$ is locally finite,
hence it is clear that a point $p$ is in the boundary of $P=\cup_{x\in\G(u)}\Qx$
if and only if $p$ is not an interior point of $\Qx$ for all $x\in\G(u)$
and $p$ is a boundary point of $\Qx$ for some $x\in\G(u)$.
\end{rem}

\begin{lem}
\label{IntFe-in-ddQu}
We have $\Int\Fe\subset\dd\Qu$.
\end{lem}

\begin{proof}
Suppose there is a point $a\in\Int\Fe=\Int(\Ee\cap\ddP)$ such that $a\not\in\dd\Qu$.
Since $a\in\dd P$ and $a\not\in\dd\Qu$ there exists $x\in\Guou$ such that $a\in\dd\Qx$.
Then any neighbourhood of $a$ intersects $\Ee\cap\Int\Qx\subset\Ee\bs\ddP$.
The projection $\covbdlmap:\tLrp\to\tGrp$ is continuous
and the restriction $\covbdlmap|_{\ddP}:\ddP\to\tGrp$ is a homeomorphism,
therefore any neighbourhood of $a$ intersects
$((\covbdlmap|_{\ddP})^{-1}\circ\covbdlmap)(\Ee\bs\ddP))\subset\ddP\bs\Ee$.
This implies $a\not\in\Int(\Ee\cap\ddP)=\Int\Fe$. Contradiction.
\end{proof}

\begin{prop}
\label{FG=Fe}
We have
$$\Fe=\Cl\Int\left((\Ee\cap\dd\Qu)-(\cupIntQxGuou)\right).$$
\end{prop}

\begin{proof}
Let $\hF:=(\Ee\cap\dd\Qu)-(\cup_{x\in\Guou}\Int\Qx)$.
It holds $\Fe\subset\Ee$ by definition, $\Int\Fe\subset\dd\Qu$ by
lemma~\ref{IntFe-in-ddQu} and $\Fe\cap\Int\Qx=\emptyset$ because of
$\Fe\subset\ddP$.
This implies $\Int\Fe\subset\hF$ and therefore $\Int\Fe\subset\Int\hF$.
On the other hand from $\hF\subset\Ee\cap\ddP$
it follows that $\Int\hF\subset\Int(\Ee\cap\ddP)=\Int\Fe$.
From $\Int\hF=\Int\Fe$ it follows that $\Fe=\Cl\Int\hF$.
\end{proof}

\begin{lem}
\label{Fe-compact}
If $\G$ is co-compact, then $\Fg$ is compact.
\end{lem}

\begin{proof}
Consider a sequence $a_k$ in~$\Int\Fg$.
Let $\varphi$ be the composition of the projection maps $\ddP\to\tGrp$ and $\tGrp\to\tGrp/\G$.
Since the quotient~$\tGrp/\G$ is compact we may assume without loss of
generality that the sequence $\varphi(a_k)$ tends to a limit $\bar a\in\tGrp/\G$.
Since $\varphi$ is surjective there exists a pre-image~$a\in\ddP$ of~$\bar a$ under~$\varphi$.
Hence there is a sequence $h_k$ in~$\G$ such that the sequence $h_k a_k$ tends to~$a$.
Since the family $(\Fg)_{g\in\G}$ is locally finite there exists a
neighbourhood~$U$ of~$a$ that intersects only finitely many fundamental domains~$\Fg$,
and therefore the set $\{h_k|k\in\n\}$ is finite.
After choosing a subsequence we may assume that the sequence $h_k$ is constant, say $h_k=h$.
Then the sequence $h a_k$ tends to~$a$,
hence the sequence $a_k$ tends to $h^{-1}a$.
This implies $h^{-1}a\in\Fg$.
\end{proof}

\begin{thmB}
If $\G$ is co-compact then  $\Fg$ is a compact polyhedron,
i.e. a finite union of finite compact intersections of half-spaces~$\Ia$.
\end{thmB}

\begin{proof}
The family $(\Qx)_{x\in\G(u)}$ is locally finite
and the fundamental domain $\Fe$ is compact by lemma~\ref{Fe-compact}.
From this it follows that there is a finite subset $E\subset\Gu$ such that
$\Fe\cap\Qx=\emptyset$ for all $x\in\Gu\bs E$.
By proposition~\ref{IntFe-in-ddQu} this implies the assertion.
\end{proof}

\section{Examples}

\label{examples}

\noindent
We have computed the fundamental domains explicitly for those infinite series of discrete subgroups,
which correspond via the construction of I.~Dolgachev~\cite{Do83}
to the Arnold series $E$, $Z$ and $Q$ of quasi-homogeneous surface singularities.
In particular the quotient of $\tSU$ by one of this groups is diffeomorphic 
to the link of the corresponding quasi-homogeneous singularity.
Whenever it is convenient, we shall denote these subgroups also by the symbols $E_n$, $Z_n$, $Q_n$.

\myskip
A discrete co-compact
subgroup $\Gamma$ of level $k$ in $\tSU$ such that the image in
$\PSU$ is a triangle group with signature
$(\alpha_1,\alpha_2,\alpha_3)$
will be denoted by~$\Gamma(\alpha_1,\alpha_2,\alpha_3)^k$.
The subgroups $E_n$, $Z_n$, $Q_n$ are of this type. 
The level $k$ and the signature $(\alpha_1,\alpha_2,\alpha_3)$ are given in the following table
(compare K.~M\"ohring~\cite{Mo}, table~19).
\begin{center}
  \setlength\XX{0.1\baselineskip}
  \setlength{\extrarowheight}{3\XX}
  \begin{tabular}{|c|c|c|c|}
  \hline
    \upshape Type &\upshape $n$~mod~$4$& $k$&$(\alpha_1,\alpha_2,\alpha_3)$
    \\[2\XX]
  \hline
           &  $0$& $(n-10)/2$& $(2,3,k+6)$\noline
      $E_n$&  $2$& $(n-10)/4$& $(3,3,k+3)$\noline
           &$1,3$& $(n-10)/3$& $(2,4,k+4)$ \\[2\XX]
  \hline
           &  $3$& $(n-9)/2$&$(2,3,2k+6)$\noline
      $Z_n$&  $1$& $(n-9)/4$&$(3,3,2k+3)$\noline
           &$0,2$& $(n-9)/3$&$(2,4,2k+4)$ \\[2\XX]
  \hline
           &  $2$& $(n-8)/2$&$(2,3,3k+6)$\noline
      $Q_n$&  $0$& $(n-8)/4$&$(3,3,3k+3)$\noline
           &$1,3$& $(n-8)/3$&$(2,4,3k+4)$ \\[2\XX]
  \hline
  \end{tabular}
\end{center}

\myskip\noindent
The following figures show some of the explicitly computed fundamental domains.
We restrict ourselves here to the series $E$.
For further figures of fundamental domains, for instance for the series $Z$ and $Q$,
for computations of fundamental domains for quadrangle rather than triangle subgroup,
and for deeper discussion of connections with quasi-homogeneous surface
singularities we refer to~\cite{BPR}, and also to~\cite{R} and \cite{Pr:2001}.

\myskip
Some explanations are required to make the figures of fundamental domains comprehensible.
The image $\pi(F_e)$ of the fundamental domain $F_e$ for a discrete co-compact
group $\Gamma\subset\tSU$ of finite level is a compact
polyhedron in~$\su(1,1)$ with flat faces.  The Lie algebra~$\su(1,1)$ is a
$3$-dimensional flat Lorentz space of signature $(n_+,n_-) = (2,1)$.  Such a
polyhedron has a distinguished rotational axis of symmetry.  The direction of
this axis is negative definite, and the orthogonal complement is positive
definite.  Changing the sign of the pseudo-metric in the direction
of the rotational axis transforms Lorentz space into a well-defined Euclidean space.
The image~$\pi(F_e)$ of the fundamental domain is then transformed into a
polyhedron in Euclidean space with dihedral symmetry.  
Figure~\ref{serE} shows the Euclidean polyhedra obtained in this way. 
The direction of the rotational axis is vertical. The top and bottom faces are removed.

The polyhedra in figure~\ref{serE} are all scaled by the same factor
to illustrate the proportions between different fundamental domains.

Figure~\ref{identE} illustrates the identification scheme 
for the cases $E_{10+2n}$.
The face identification is equivariant with respect to
the dihedral symmetry of the polyhedron.  The faces shaded in the same way are
identified.  Arrows on the edges of shaded faces indicate the identified flags
(face, edge, vertex).

\def\po #1#2{\off{po#2_#1}}

\begin{figure}
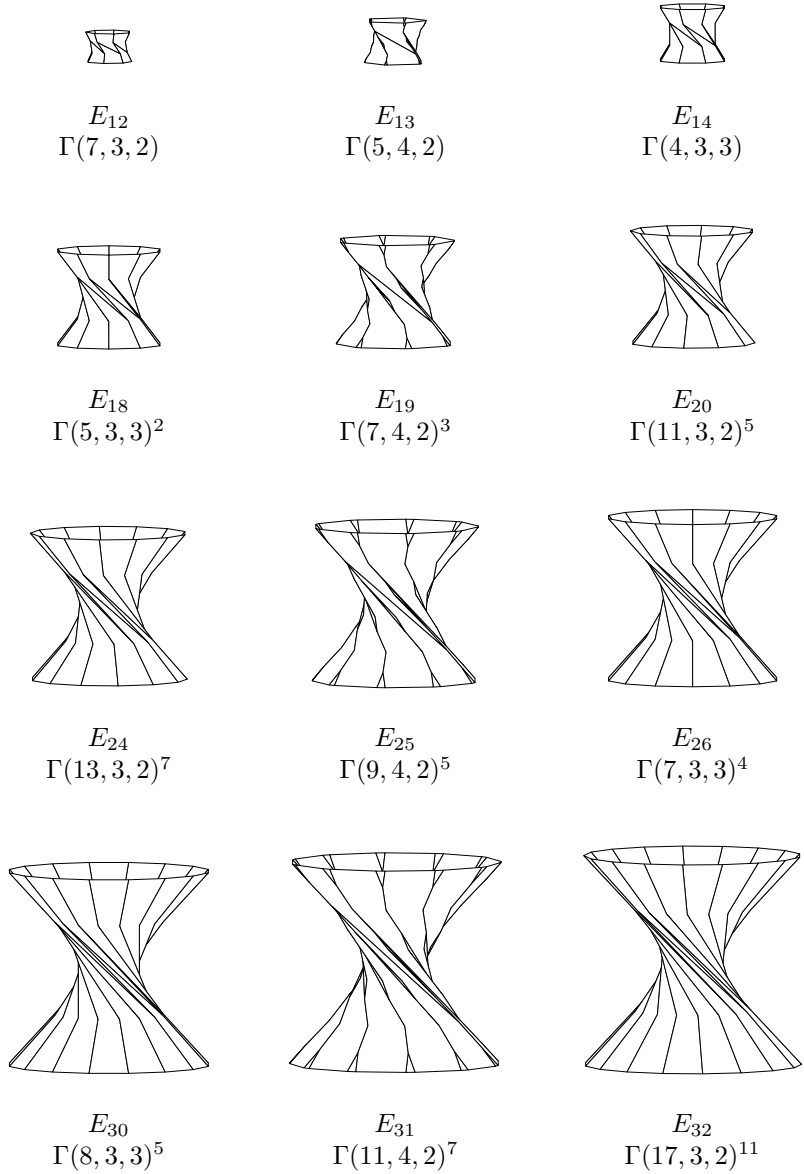

\begin{center} \def\off#1{\epsffile{t5.#1.ps}}
\begin{tabular}{ccc}
    \po E{237}         & \po E{245}        & \po E{334} \\
    $E_{12}$           & $E_{13}$          & $E_{14}$   \\
    $\Gamma(7,3,2)$    & $\Gamma(5,4,2)$   & $\Gamma(4,3,3)$ \\\\
    \po E{335}         & \po E{247}        & \po E{2311} \\
    $E_{18}$           & $E_{19}$          & $E_{20}$   \\
    $\Gamma(5,3,3)^2$  & $\Gamma(7,4,2)^3$ & $\Gamma(11,3,2)^5$ \\\\
    \po E{2313}        & \po E{249}        & \po E{337} \\
    $E_{24}$           & $E_{25}$          & $E_{26}$   \\
    $\Gamma(13,3,2)^7$ & $\Gamma(9,4,2)^5$ & $\Gamma(7,3,3)^4$ \\\\
    \po E{338}         & \po E{2411}       & \po E{2317} \\
    $E_{30}$           & $E_{31}$          & $E_{32}$   \\
    $\Gamma(8,3,3)^5$  & $\Gamma(11,4,2)^7$& $\Gamma(17,3,2)^{11}$
\end{tabular}
\end{center}
\caption{Fundamental domains for the beginning of the $E$-series}
  \label{serE}
\end{figure}

\input identE.pic

\section{Concluding Remarks}

\label{outlook}

\begin{enumerate}[1)]
\item
The construction of fundamental domains in the flat Lorentz case using crooked planes
by T.~Drumm and W.~Goldman 
(see~\cite{DG95}, \cite{DG99} and references therein)
is the only other fundamental domain construction in the non-Rie\-mann\-ian pseudo-Riemannian case we are aware of,
besides the construction for $\tSU$ presented in this paper
and the construction for $\PSU$ studied before in \cite{Fi}, \cite{KNRS}, \cite{Ba}.
\item
The idea of projection of an affine construction with half-planes onto a quadric
is also used in the algorithmic construction of Voronoi diagrams for (finite) point sets
in the Euclidean and hyperbolic plane, compare J.-D.~Boissonnat and M.~Yvinec~\cite{BY}, Part~V.
\item
R.S.~Kulkarni and F.~Raymond studied in~\cite{KR} standard Lorentz space forms.
A Lorentz space form is called standard if it is a quotient of $\tGrp=\tSU$ 
by a discrete subgroup of $\Isom(\tGrp)$ conjugate to a subgroup of $J=\<J_0,\eta\>$,
where $J_0=(G\times K)/\De_Z\subset\Isom_0(\tGrp)$ and $K=\{r_0(t)\}_{t\in\r}$.
We recall that $\Isom(\tGrp)=\<\Isom_0(\tGrp),\spieg,\inv\>$ 
and $\Isom_0(\tGrp)\cong(\tGrp\times\tGrp)/\De_Z$.
Examples of non-standard Lorentz space forms were found by W.~Goldman~\cite{Go}
and recently by F.~Salein~\cite{Sa}.

\myskip
We can think of a discrete subgroup of $\tGrp$ acting by left translations 
as a discrete subgroup of $(G\times\{e\})/\De_Z\subset J_0$,
so the Lorentz space forms studied in this paper are standard.
We would like to generalize our fundamental domain construction
for the case of other Lorentz space forms, at least for standard ones.

\myskip
For instance, our fundamental domain construction can be modified
to work for the subgroups of $\Isom_0(\tGrp)$ of the form $(\G\times\Phi)/\De_Z$,
where $\G$ is a discrete subgroup of~$\tGrp$ and $\Phi$ is a discrete subgroup of $K$.
Also the work of I.~Dolgachev~\cite{Do83} can be generalized to obtain
a correspondence between subgroups of this special form and quasi-homogeneous
$\q$-Gorenstein surface singularities.
These constructions will be studied in a forthcoming paper.
\item 
One can also think of generalizations of the fundamental domain construction
in the cases of other (pseudo-)Riemannian quadrics in pseudo-Euclidean spaces. 
In section~\ref{analogons} we discussed the analogous constructions for two one-dimensional cases. 
It would be interesting to generalize the constructions
for higher dimensional quadrics.
\end{enumerate}

\bibliographystyle{amsalpha}
\bibliography{statia}

\end{document}